\documentclass[journal]{IEEEtran}
\IEEEoverridecommandlockouts 

\usepackage{cite}
\usepackage{amsmath, amstext, amssymb}
\usepackage{array}
\usepackage{stfloats}
\usepackage{url}

\ifCLASSINFOpdf
   \usepackage[pdftex]{graphicx}
   \graphicspath{{./Images/}}
   \DeclareGraphicsExtensions{.pdf,.jpeg,.png}
\else
\fi

\usepackage{multirow}
\usepackage{booktabs,colortbl}

\begin{document}
%
\title{Corrective Control to Handle Forecast Uncertainty: \\ A Chance Constrained Optimal Power Flow}

\author{Line~Roald,~
        Sidhant Misra,~
        Thilo Krause,~
        and~G\"oran~Andersson~


%

\thanks{The work described in this paper has been partially funded by the project "Innovative tools for future coordinated and stable operation of the pan-European electricity transmission system (UMBRELLA)", supported under FP7 of the European Union, grant agreement 282775. }
\thanks{L. Roald, T. Krause and G. Andersson are with the Power System Laboratory, ETH Zurich, Switzerland. E-mail: roald@eeh.ee.ethz.ch}
\thanks{S. Misra is with the Center for Nonlinear Studies and Theoretical Division T-4 of Los Alamos National Laboratory, Los Alamos, United States.}}


\maketitle

\begin{abstract}
Higher shares of electricity generation from renewable energy sources and market liberalization is increasing uncertainty in power systems operation. At the same time, operation is becoming more flexible with improved control systems and new technology such as phase shifting transformers (PSTs) and high voltage direct current connections (HVDC). 
Previous studies have shown that the use of corrective control in response to outages contributes to a reduction in operating cost, while maintaining N-1 security. In this work, we propose a method to extend the use of corrective control of PSTs and HVDCs to react to uncertainty. We characterize the uncertainty as continuous random variables, and define the corrective control actions through affine control policies. This allows us to efficiently model control reactions to a large number of uncertainty sources. 
The control policies are then included in a chance constrained optimal power flow formulation, which guarantees that the system constraints are enforced with a desired probability. By applying an analytical reformulation of the chance constraints, we obtain a second-order cone problem for which we develop an efficient solution algorithm. 
In a case study for the IEEE 118 bus system, we show that corrective control for uncertainty leads to a decrease in operational cost, while maintaining system security. Further, we demonstrate the scalability of the method by solving the problem for the IEEE 300 bus and the Polish system test cases.
\end{abstract}

\begin{IEEEkeywords}
Corrective Control, Renewable Integration, Chance Constraints, Optimal Power Flow
\end{IEEEkeywords}

\section{Introduction}
In power systems terminology, a \emph{preventively} secure system state refers to a situation in which the power system remains secure after any credible contingency (typically defined as an N-1 situation) without any additional control action. A \emph{correctively} secure state refers to a situation where additional, post-disturbance controls might be required.

The use of corrective security has increased over the past years \cite{panciatici2014}. On the one hand, the use of corrective actions is driven by increased uncertainty from renewable electricity generation, economic considerations due to market liberalization and operation of the system closer to its limits. On the other hand, better control systems, more situational awareness and the installation of devices such as phase-shifting transformers (PSTs) and high-voltage direct current (HVDC) connections provide the system operator with new possibilities to control power flows and react to changes in the system in real-time. 


While corrective control can reduce the operational cost \cite{chatzivasileiadis2011} and is applied routinely 
in real-time system operation, 
the real-time set-point changes of HVDC and PSTs are typically chosen in an ad-hoc fashion. This is partially due to the difficulty of planning corrective control actions in response to forecast uncertainty, as it requires the consideration of a large number possible uncertainty scenarios in addition to consideration of contingencies. However, the need to ensure that corrective control will be sufficient in real-time calls for efficient ways of modelling both the possible system states and the corresponding corrective control actions.

Power system operational planning with consideration of uncertainty has been considered in many different ways, e.g. 
\cite{bouffard2008, sjodin2012, vrakopoulou2012, roald2013, bienstock2014, warrington2013, lorca2015, li2015, summers2014}, but only few have considered the application of corrective control actions or the existence of power flow control devices such as PSTs or HVDC. 
In \cite{panciatici2010, capitanescu2012}, a three-stage optimal power flow (OPF) framework where corrective control actions are used to ensure feasibility during worst-case combinations of contingencies and uncertainty was proposed, and was extended to the definition of optimal corrective control actions in \cite{fliscounakis2013}. 
The OPF formulation in \cite{mueller2014} accounts for uncertainty and includes corrective control for a limited number of pre-selected, critical uncertainty scenarios. 
In \cite{vrakopoulou2013ISGT, thakurta2015}, \emph{post-contingency} corrective control of power flow control devices such as HVDC and PSTs is applied within a chance constrained OPF (CC-OPF) framework, which ensures that the constraints will hold with a desired probability. However, corrective control actions \emph{in reaction to the uncertainty realizations} have not been considered in any of these works. 

In this paper, we extend previously proposed CC-OPFs to account for corrective control in reaction to \emph{both contingencies and uncertainty}.
Corrective control for uncertainty differs from post-contingency corrective control in several ways. While contingencies are typically low probability, discrete events that induce large and sudden impacts on the system, uncertainty realizations occur frequently and develop in a more continuous fashion (although ramping due to, e.g., sunset or fog can be fast). 
These differences impact the time available for implementation, as well as the type and modelling of control reactions. 
However, we refer to both as \emph{corrective control}, since they act to mitigate the impact of already occurred events (as opposed to preventive measures). 

We focus on corrective control of HVDC and PSTs, which are typically controlled by the transmission system operator and have low cost, as opposed to the use of generation redispatch, which interferes with market operation and incurs significant cost.
Further, we choose to work with analytically reformulated chance constraints to account for the impact of uncertainty, as these offer a transparent and scalable way of ensuring security with a large number of uncertainty sources. Chance constraints also align well with several methods applied in industry, such as the probabilistic reserve dimensioning applied in ENTSO-E 
\cite{ENTSOE2013supportingLFCR} or the definitions of reliability margins for flow-based market coupling in parts of Europe \cite{CWE2011}.
The contributions of the paper can be summarized in the following points:
\begin{enumerate}
    \item We propose a framework with combined corrective control in reaction to both uncertainty and contingencies. Since uncertainty from, e.g., renewable energy production are naturally characterized as continuous random variables, we propose to model the corrective control through a continuous, affine control policy. 
    \item We formulate a chance constrained optimal power flow with security constraints (CC-SCOPF) based on a combination of the formulations developed in \cite{roald2013, bienstock2014}, and extend it to include the proposed corrective control framework.
    \item The CC-SCOPF is reformulated into an optimization problem with second-order cone (SOC) constraints, for which we develop an efficient solution algorithm. The proposed algorithm is based on solving a sequence of second-order cone programs (SOCPs), which in our case outperforms the cutting-plane algorithm proposed in \cite{bienstock2014}.
\end{enumerate}
The benefits of the proposed CC-SCOPF are demonstrated in a case study based on the IEEE 118 bus system. The results show that corrective control in reaction to uncertainty can reduce operational cost, while maintaining system security. Further, we demonstrate and discuss the scalability of the method in a case study including both the IEEE 300 bus test system and the large-scale Polish test case with 2383 buses. 

The remainder of the paper is organized as follows. Section \ref{sec:corrcontrol} describes the general framework of corrective control, based on a generic power flow controller. The power system modelling, including uncertainty and corrective control, is described in Section \ref{sec:modelling}. Section \ref{sec:CCOPF} provides the complete CC-SCOPF formulation and discusses the reformulation of the chance constraints. Details of the sequential SOCP algorithm are given in Section \ref{sec:implementation}. 
The case study for the IEEE 118 bus system is included in Section \ref{sec:Case}, 
while Section \ref{sec:Polish} demonstrates and discusses scalability of the method on larger test cases. Section \ref{sec:Conclusion} summarizes and concludes the paper.

\section{Modelling Framework for Corrective Control}
\label{sec:corrcontrol}
In this paper, we consider corrective control as a means to handle transmission line congestion. We distinguish between two types of corrective control: Corrective control to handle contingencies and corrective control to handle forecast uncertainty. 
In the following, we present how we model corrective control for a generic power flow controller, which changes set-point in response to transmission line outages or fluctuations in the power injections. 
Corrective control for generation outages is not considered here, but can be incorporated within the proposed modelling framework.
Modelling considerations related specifically to HVDC and PSTs, as well as how the corrective control influences the power flows on the transmission lines, are given in Section \ref{sec:modelling}.

We denote vectors by lower case letters, e.g., $p_G,~\omega$. The components of the vectors are denoted by using lower case subscripts, i.e, the $i$th component of $p_G$ is denoted by $p_{G,i}$.
Matrices are denoted by upper/lower bold case letters, $\boldsymbol{\alpha}, \mathbf{M}$, and $\boldsymbol{\alpha}_{(i,\cdot)}, \boldsymbol{\alpha}_{(\cdot,i)}$ denote the $i^{th}$ row and column of $\boldsymbol{\alpha}$ respectively.
Index $i$ refers to generators, while index $ij$ refers to lines.

\subsection{Modelling Corrective Control for Contingencies}
Post-contingency control is modeled as a change in the set-points of power flow control devices (such as the flow across HVDC connection or the phase angles of PSTs) following an outage. 
Since outages can be characterized as a set of discrete events, we model the set-point changes by introducing additional optimization variables $\delta^{ij}$ in for every outage as in \cite{chatzivasileiadis2011, vrakopoulou2013ISGT}.
The set-point of a power flow controller following the outage of line $ij$ is thus given by 
\begin{equation}
\tilde{\pi}(ij) = \pi + \delta^{ij},
\label{policy}
\end{equation}
where $\pi,~\tilde\pi(ij)$ represents the pre- and post-contingency set-points, respectively, and $\delta^{ij}$ represents the set-point change following the contingency. 

\subsection{Modelling Corrective Control for Forecast Uncertainty}
While post-contingency corrective control describes reactions of the controllers to an outage \emph{after the outage happens}, corrective control for uncertainty describes reactions of the controllers \emph{after the uncertainty has realized}. 
We model these corrective control actions through affine control policies, where the controllers adjust their set-point proportional to the forecast error. The controller set-points, including corrective control for forecast uncertainty, are given by 
\begin{equation}
\tilde{\pi}(\omega) = \pi + \boldsymbol{\alpha} \omega,
\label{policyUnc}
\end{equation}
where $\pi,~\tilde\pi(\omega)$ represents the scheduled and real-time controller set-points, and the vector $\omega$ denotes the random fluctuations. The matrix  $\boldsymbol{\alpha}$ defines the response of the controllers, and control parameters $\boldsymbol{\alpha}$ are subject to optimization. 

Restricting the corrective action to follow a linear response leads to solutions that have higher cost than if we use, e.g., an arbitrary function or define optimal set-points based on specific uncertainty scenarios. However, representing arbitrary functions or a reasonable set of uncertainty scenarios within an optimization problem is challenging, and negatively affects scalability to large systems. 
Therefore, while an affine response policy 
restricts the ability to react, 
it also has several advantages. First, it allows us to optimize over the corrective control actions without compromising computational tractability.
Second, a policy-based response allows us to treat the fluctuations as continuous variables (as opposed to a representation through a finite number of scenarios). Third, it provides a control policy which can be easily implemented by the TSO. 
Fourth, even if the policy is not directly implemented, it guarantees that a feasible solution will exist. 

Affine policies have been utilized within chance constrained and robust models for the optimal power flow problem with uncertainty, e.g. \cite{vrakopoulou2012, roald2013, warrington2013, bienstock2014, lorca2015}, to model the continuous activation of reserves through a mechanism similar to the automatic generation control (AGC).
The approach taken in this paper differs from previous work in two important aspects. 
First, the affine control policies defined previously were used to \emph{balance} the grid (by controlling generation), while we focus on corrective power flow control for \emph{congestion management} (by controlling HVDC and PSTs). 
Second, previously proposed policies focused on balancing the overall power mismatch without consideration of the location of the fluctuations in the grid.
Since we apply corrective control to manage power flows and congestion, the specific location of both the fluctuations and the control reactions in the grid matters. We therefore define control parameters $\boldsymbol{\alpha}_{(i,j)}$ for each pair of controller $i$ and uncertainty source $j$, such that the controllers can react to fluctuation $\omega_j$ based on its location in the grid.

\subsection{Combined Corrective Control}
Assuming that each controller is able to provide corrective control for both contingencies and forecast uncertainty, we obtain the combined corrective control policy,
\begin{equation}
\tilde{\pi}(\omega,ij) = \pi + \boldsymbol{\alpha} \omega + \delta^{ij}.
\label{policy_combined}
\end{equation}
Note that the reaction to forecast uncertainty is decoupled from the reaction to contingencies, i.e., $\boldsymbol{\alpha}$ does not depend on whether the system is in post- or pre-contingency state. 
This is because the number of additional control variables would otherwise be very significant and difficult to handle, both for the optimization and for implementation in real-time operation. 

\subsection{Use of Corrective Control in Real-Time Operation}
The corrective control actions $\delta^{ij},\boldsymbol{\alpha}$ are necessary to ensure secure system operation in real time. 
The post-contingency corrective control actions $\delta^{ij}$ can be made available to the operator in form of a look-up table, such that they can be implemented fast following a contingency. 
The corrective control actions due to uncertainty, defined by a combination of the real time realization of the fluctuations $\omega$ and the control parameters $\boldsymbol{\alpha}$, can be implemented 
either as continuous, automatic adjustments (e.g., similar to AGC control)
or as a manual change in set-points which is implemented only when a constraint violation occurs in the system. 

While restricting the response to an affine policy allows for a tractable day-ahead OPF, we can implement any reaction in real-time operation. In real time the uncertainty $\omega$ is known, and the operator can rerun the SCOPF as a deterministic problem. The set-points of the HVDC and PSTs can then be optimally chosen based on the current operating point. Since the real-time control is more general than the affine policy, we solve an SCOPF with less restrictive constraints, extending the feasible space. In this regard, applying an affine policy in operational planning guarantees feasibility in real time operation.


\section{Power System Modelling} \label{sec:modelling}
In this section, we present details about the power system modelling. We consider a power system where $\mathcal{N},~\mathcal{L}$ is the set of nodes and lines, respectively, and the number of nodes and lines are given by $|\mathcal{N}| = m$ and $|\mathcal{L}| = l$.
The set of nodes with uncertain demand or production of energy is given by $\mathcal{U}\subseteq\mathcal{N}$. The fluctuations in the demand or production at any given node can be due to various sources, such as load fluctuations, forecast errors for wind or PV or intra-day electricity trading.  
The set of conventional generators is denoted by $\mathcal{G} \subseteq \mathcal{N}$, and are assumed to be controllable within their limits.
To simplify notation, we assume that there is one conventional generator $p_{G,i}$, one composite uncertainty source  $u_i$ and one demand $d_i$ per node, such that $|\mathcal{G}|=|\mathcal{U}|=|\mathcal{N}| = m$. Nodes without generation or load can be handled by setting the respective entries to zero, and nodes with multiple entries can be handled through a summation.
The considered power flow control devices used for corrective control are given by the set of HVDC connections $\mathcal{H}$ and PSTs $\mathcal{S}$, with $|\mathcal{H}| = h$ and $|\mathcal{S}| = s$, respectively.
The modelling is based on a DC power flow approximation. We consider the outage of any single line for the N-1 security constraints, but leave out lines that lead to islanding of the system. Generation outages are handled in a simplified way through a pre-determined reserve requirement.

\subsection{Uncertainty Sources}
The uncertainty sources $u\in\mathbb{R}^m$ are modeled as the sum of the expected production of active power (e.g., from wind or solar PV) $\mu$ and a zero mean fluctuating component $\omega$:
\begin{equation}
u = \mu + \omega~.
\end{equation}
The covariance matrix of the fluctuations is denoted by $\Sigma_W \in \mathbb{R}^{m\times m}$. Further, the total power mismatch arising from the fluctuations is given by $\Omega=\sum_{i\in\mathcal{U}}\omega_i$, with corresponding standard deviation $\sigma_{\Omega} = \sqrt{\mathbf{1}^{T}\Sigma_W \mathbf{1}}$.

\subsection{HVDC and PST Set-Points With Corrective Control}
We assume that all HVDC installations are point-to-point connections
such that the power flow $p_{DC}$ is controllable within the limits  $\underline{p}_{DC}$, $\overline{p}_{DC}$. 
The PST tap positions are assumed to be close enough for the angle $\gamma$ to be well approximated as continuous variables, and the lower and upper limit on the PST angles are given by $\underline{\gamma}, \overline{\gamma}$.
As outlined in Section \ref{sec:corrcontrol}, the post-contingency corrective control of HVDC and PSTs are modelled through additional variables $\delta_{DC}^{ij},~\delta_{\gamma}^{ij}$, while the corrective control for uncertainties are modelled through the matrices $\boldsymbol{\alpha}_{DC}\in\mathbb{R}^{h \times m},~\boldsymbol{\alpha}_{\gamma}\in\mathbb{R}^{s \times m}$, respectively. The set-points of the HVDC and the PSTs with corrective control are thus given by
\begin{align}
&\tilde{p}_{DC,i}(\omega,ij)=p_{DC,i} + \delta_{DC,i}^{ij} - \boldsymbol{\alpha}_{DC(i,\cdot)}\omega, \quad \forall_{i\in\mathcal{H}, ij\in \{0,\mathcal{L}\}}, \label{HVDCcontrol}\\
&\tilde{\gamma}_{DC,i}(\omega,ij)=\gamma_{i} + \delta_{\gamma,i}^{ij} - \boldsymbol{\alpha}_{\gamma(i,\cdot)}\omega, \quad \forall_{i\in\mathcal{S}, ij\in \{0,\mathcal{L}\}}, \label{PSTcontrol}
\end{align}
where $ij=0$ refers to the set-points in the pre-contingency state where $\delta_{DC,i}^{0}=\delta_{\gamma,i}^{0}=0$.

\subsection{Power Balance and Generation Control}

One of the main tasks in power system operation is to ensure balance between consumed and produced power at all times. Here, we split the modelling of the power balance into the case with and without deviations $\omega$ from the schedule. For the base case (with $\omega=0$ and no outages), we enforce the nodal power balance constraints,
\begin{equation}
p_G - d + \mu + \mathbf{C}_{DC}p_{DC}= \mathbf{B}_{Bus}\cdot\theta + \mathbf{B}_{\gamma}\cdot\gamma~. \label{nodal_bal_first}
\end{equation} 
Here, $\mathbf{B}_{Bus}\in\mathbb{R}^{m\times m}$ is the bus susceptance matrix of the system and $\theta\in\mathbb{R}^m$ refers to the voltage angles at each bus. 
$\mathbf{C}_{DC}$ is the incidence matrix of the HVDC connections, with -1 at the bus where the connection is leaving and +1 at the bus where it enters.
The matrix $\mathbf{B}_{\gamma}$ describes the influence of the PSTs and has non-zero entries only for buses connected to lines with PSTs. For PST $k$ located at line $ij$ (which leaves from bus $i$ and enters bus $j$), we have the following non-zero elements, 
\begin{equation}
\mathbf{B}_{\gamma(i,k)}=\frac{1}{x_{ij}},\quad \mathbf{B}_{\gamma(j,k)}=-\frac{1}{x_{ij}}.
\end{equation} 

To balance fluctuations $\omega$, we assume that the generators adjust their in-feeds according to an automatic generation control (AGC) signal. 
We base the activation of reserves on the total power mismatch $\Omega$ as in \cite{vrakopoulou2012, roald2013, warrington2013, bienstock2014, lorca2015}. In this case, the generator set-points after adjustment for the fluctuations are given by 
\begin{equation}
\tilde{p}_G(\omega) = p_G - \alpha_G \Omega, \label{global}
\end{equation}
where $\alpha_G \in \mathbb{R}^m$ is a vector distributing the power mismatch among the generators, and is subject to optimization. 
To ensure active power balance during fluctuations and that the reserve activation contributes to balancing of the deviation (i.e., a positive $\Omega$ induce a decrease in other generation), 
we additionally enforce the following constraints on $\alpha_G$:
\begin{equation}
\sum_{i\in\mathcal{G}} \alpha_{G,i} = 1, \quad \alpha_G \geq 0. \label{windbal}
\end{equation}

Note that we do not include corrective control by generators for neither contingencies nor uncertainties. While the formulation itself can easily be extended by variables denoting post-contingency redispatch $d_G^{ij}$ or a more general matrix $\boldsymbol{\alpha}_G$, generation control is expensive and require careful consideration of how redispatch and balancing is priced. Allowing a more general matrix $\boldsymbol{\alpha}_G$ can however significantly impact the handling of uncertainty as shown in \cite{roald2016EnergyCon}.\\
Since line outages do not change the power injections and power balance is maintained during fluctuations due to \eqref{windbal}, the above modelling enforces power balance during all considered system conditions.

\subsection{Power flow modelling}
We compute the pre-contingency power flows as the sum of the base case flow $p_{ij}$ (with $\omega=0$ and no contingencies) and changes due to fluctuations $\Delta p_{ij}^{\omega}$: 
\begin{equation}
p_{ij}^{\omega} = p_{ij} + \Delta p_{ij}^{\omega}, 
\end{equation}
The base case power flow $p_{ij}$ is given by 
\begin{equation}
p_{ij} = \frac{1}{x_{ij}}( \theta_i - \theta_j) + b_{\gamma(ij,\cdot)}\gamma, 
\end{equation}
where $x_{ij}$ is the reactance of lines $ij$. The matrix $b_{\gamma}\in\mathbb{R}^{l\times s}$ maps the PST angles to the lines where the PSTs are located, with entries $b_{\gamma(ij,s)}=\frac{1}{x_{ij}}$ if the $s^{th}$ PST is located at line $ij$ and $b_{\gamma(ij,s)}=0$ otherwise. 

The power flow change $\Delta p_{ij}^\omega$ is a result of both the fluctuations $\omega$ themselves, the corrective control from HVDC and PSTs and the balancing by the generators. It is given by
\begin{align}
& \Delta p_{ij}^{\omega} = \mathbf{M}_{(ij,\cdot)}\left(-\alpha_G 1_{1\times m}  + \mathbf{I} + \mathbf{C}_{DC} \boldsymbol{\alpha}_{DC} - \mathbf{B}_{\gamma} \boldsymbol{\alpha}_{\gamma}\right)\omega~  \nonumber \\
& \quad \quad \quad ~ +  b_{\gamma(ij,\cdot)}\boldsymbol{\alpha}_{\gamma}\omega  = \mathbf{A}_{(ij,\cdot)}\omega~,  \label{eq:lineflows1}
\end{align}
where $\mathbf{I}\in \mathbb{R}^{m \times m}$ is the identity matrix. 
The matrix  $\mathbf{M}\in\mathbb{R}^{l \times m}$ is the matrix of power transfer distribution factors (PTDFs) \cite{christie2000}.
The term in the brackets is the effective change in the power injections due to fluctuations $\omega$, and the last term is the direct influence of the changes in the PST set-points on the lines where they are located.
Note that the matrix $\mathbf{A} = \mathbf{A}(\alpha_G, \boldsymbol{\alpha}_{DC}, \boldsymbol{\alpha}_{\gamma})$ is a linear function of the policies for balancing and corrective control. By allowing corrective control of HVDC and PSTs in reaction to uncertainties, it is possible to influence the change in the power flow and reduce line congestion.

With line outages, the power flow changes both due to corrective control actions from HVDC and PSTs, and due to the change in system topology. 
The change in the power flow due to the corrective control actions alone (including the change that would have occured on the outaged line itself) can be modelled as
\begin{align}
p_{ij}^{\delta_{kl}} &= p_{ij} + \mathbf{M}_{(ij,\cdot)}\mathbf{C}_{DC} \delta_{DC}^{kl} + \left(b_{\gamma(ij,\cdot)} - \mathbf{M}_{(ij,\cdot)}\mathbf{B}_{\gamma}\right) \delta_{\gamma}^{kl}. \nonumber
\end{align}
The effect of change in system topology can be accounted for using line outage distribution factors (LODFs) \cite{christie2000}, with $LF_{ij}^{kl}$ denoting the  fraction of the power flow on the line $kl$ that is shifted to line $ij$ when line $kl$ is outaged.
The flow on a line $ij$ with outage $kl$ and fluctuation $\omega$ is given by 
\begin{align}
p_{ij}^{kl,\omega} &= p_{ij}^{\delta_{kl}} + \Delta p_{ij}^{\omega} + LF_{ij}^{kl}\left( p_{kl}^{\delta_{kl}} + \Delta p_{kl}^{\omega}\right)  \nonumber\\
&= p_{ij}^{\delta_{kl}} + \mathbf{A}_{(ij,\cdot)}\omega + LF_{ij}^{kl}\left( p_{kl}^{\delta_{kl}} + \mathbf{A}_{(kl,\cdot)}\omega \right)  \nonumber\\ 
&= p_{ij}^{\delta_{kl}} + LF_{ij}^{kl} p_{kl}^{\delta_{kl}} + \left(\mathbf{A}_{(ij,\cdot)}+LF_{ij}^{kl}\mathbf{A}_{(kl,\cdot)}\right) \omega. \nonumber 
\end{align}


\section{Chance-Constrained Optimal Power Flow} \label{sec:CCOPF}
In this section, we formulate the chance constrained optimal power flow problem based on the modelling considerations described in Section \ref{sec:modelling}.

\subsection{Objective and Constraints}

\subsubsection{Objective}
The objective of the CC-OPF is to minimize the total cost of energy and reserves,
\begin{equation} \label{eq:objective}
\min \sum_{i \in \mathcal{G}}\left(c_i p_{G,i} + c_i^+ r_i^+ + c_i^- r_i^-\right)
\end{equation}
Here, $p_{G,i}$ are the scheduled generation set-points, $r^+, r^-$ represent the up- and down regulation reserves and the cost coefficients $c_i, c_i^+, c_i^-$ represent the bids of the generators for providing energy and reserves.

\subsubsection{Power Balance and Generator Constraints}
The power balance and generator constraints are given by
\begin{align}
& p_G - d + \mu + \mathbf{C}_{DC}p_{DC}= \mathbf{B}_{Bus}\cdot\theta + \mathbf{B}_{\gamma}\cdot\gamma~, \label{nodalbal}\\
& \sum_{i\in\mathcal{G}} \alpha_{G,i} = 1,  \quad \alpha_G \geq 0, \label{windbal1}\\
& p_G + r^+ \leq \bar{p}_G, \quad  p_G - r^- \geq \underline{p}_G, \label{gen_down}\\
& \sum_{i\in\mathcal{G}} r_i^+ \geq R^+, ~~\sum_{i\in\mathcal{G}} r_i^- \geq R^-, \label{res_cont1}\\
& 0\leq r^+ \leq \bar{r}^+,~~ 0\leq r^-\leq \bar{r}^-, \label{res_cont}\\
& \mathbb{P}\left[- \alpha_{G,i}\Omega \leq r^{+}_i \right] \geq 1-\epsilon_G,~~ \forall_{i \in \mathcal{G}} \label{chance_gen_min} \\
& \mathbb{P}\left[- \alpha_{G,i}\Omega \geq r^{-}_i \right] \geq 1-\epsilon_G,~~ \forall_{i \in \mathcal{G}} \label{chance_gen_max}
\end{align}
Here, \eqref{nodalbal} defines the power balance in normal operation, while \eqref{windbal1} ensures active power balance during wind power fluctuations and non-negativity of $\alpha_G$.
The generator constraints \eqref{gen_down} ensure that the initial generation set-points in combination with scheduled reserve capacities $r^+,r^-$ remain within the technical generation limits $\bar{p}_G,~\underline{p}_G$.
Eq. \eqref{res_cont1} enforce that the total amount of reserves fulfill a predetermined reserve requirement $R^+, R^-$. Eq. \eqref{res_cont} ensure that the generators are not scheduled to provide more reserves than their ramping capabilities allow, by enforcing an upper bound on the reserve capacities $\bar{r}^+,~\bar{r}^-$ for each generator. 
Eqs. \eqref{chance_gen_min}, \eqref{chance_gen_max} enforce that the reserve activation requested from each generator 
can be covered by the corresponding reserves $r^+,~r^-$. Since the reserve activation depends on the fluctuations $\Omega$, these constraints are formulated as chance constraints, which require the probability of constraint violation to remain below an acceptable level $\epsilon_G$. \\
Note that all reserve chance constraints depend only on the total mismatch $\Omega$, which is a scalar random variable. Assuming that we enforce all generator constraints with the same acceptable violation probability $\epsilon_G$, the reserve constraints \eqref{chance_gen_min}, \eqref{chance_gen_max} are jointly enforced, i.e., the probability that \emph{none} of the generation constraints will be violated is $1-\epsilon_G$. 
To see this, we observe that the reserves of all generators will be fully utilized when the $1-\epsilon_g$ quantile $q_{1-\epsilon_g}$ of $\Omega$ is reached. For $\Omega<\hat{\Omega}_{1-\epsilon_G}$, \emph{none} of the constraints \eqref{chance_gen_max} are violated. For $\Omega>\hat{\Omega}_{1-\epsilon_G}$, \emph{all} constraints \eqref{chance_gen_max} are violated. Using this observation we have
\begin{align*}
    &\mathbb{P}\left[- \alpha_{G,i}\Omega \leq r^{+}_i \right] = \mathbb{P}\left[- \Omega \leq \frac{r^{+}_i}{\alpha_{G,i}} \right] \geq 1-\epsilon_G.
\end{align*}
Since there is a non-zero cost $c_i^{+}$ associated with $r_i^{+}$ in \eqref{eq:objective}, at optimality we will have
\begin{align*}
    \frac{r^{+}_i}{\alpha_{G,i}} = q_{1-\epsilon_g}~~ \forall_{i \in \mathcal{G}}.
\end{align*}
Hence the \emph{joint} probability of reserve insufficiency can be simplified as
\begin{align*}
    \mathbb{P}\left[- \alpha_{G,i}\Omega \leq r^{+}_i, ~~ \forall_{i \in \mathcal{G}} \right] 
    & = \mathbb{P}\left[ -\Omega \leq q_{1-\epsilon_g} \right]  = 1-\epsilon_G.
\end{align*}
%
The violation probability $\epsilon_g$ thus has a direct interpretation as the probability of having insufficient reserves, which is a commonly used risk measure in real power systems. In the Swiss power system, for example, the acceptable probability of reserve deficiency is 0.2\% \cite{abbaspourtorbati2016}, which corresponds to $\epsilon_g = 0.001$ for up and down reserves.

\subsubsection{HVDC and PST Constraints}
The constraints enforcing the upper and lower bounds on HVDC and PSTs, after activation of corrective control, are given by \eqref{chance_dc_min} -  \eqref{chance_pst_max}, while \eqref{max_corr} allows the operator to limit the amount of post-contingency control:
\begin{align}
& \mathbb{P}\left[p_{DC,i} + \delta_{DC,i}^{kl} - \boldsymbol{\alpha}_{DC(i,\cdot)}\omega \leq \bar{p}_{DC,i} \right] \geq 1-\epsilon,  \label{chance_dc_min} \\
& \mathbb{P}\left[p_{DC,i} + \delta_{DC,i}^{kl} - \boldsymbol{\alpha}_{DC(i,\cdot)}\omega \geq \underline{p}_{DC,i} \right] \geq 1-\epsilon,  \label{chance_dc_max} \\
& \mathbb{P}\left[\gamma_{j} + \delta_{\gamma,j}^{kl} - \boldsymbol{\alpha}_{\gamma(j,\cdot)}\omega \leq \bar{\gamma}_{j} \right] \geq 1-\epsilon, \label{chance_pst_min} \\
& \mathbb{P}\left[\gamma_{j} + \delta_{\gamma,j}^{kl} - \boldsymbol{\alpha}_{\gamma(j,\cdot)}\omega \geq \underline{\gamma}_{j} \right] \geq 1-\epsilon, \label{chance_pst_max} \\
& -\overline{\delta}_{DC,i}^{kl} \leq  \delta_{DC,i}^{kl} \leq -\overline{\delta}_{DC,i}^{kl},\quad 
 -\overline{\delta}_{\gamma,i}^{kl} \leq \delta_{\gamma,i}^{kl} \leq -\overline{\delta}_{\gamma,i}^{kl} \label{max_corr}\\
& \quad\quad\quad\quad\quad\quad\quad\quad\quad\quad\quad\quad\quad  \forall_{i \in \mathcal{H}, ~j \in \mathcal{S},~ kl \in \{0,\mathcal{L}\}}~. \nonumber
\end{align}
The constraints \eqref{chance_dc_min} - \eqref{chance_pst_max} depend on the uncertainty $\omega$, and are thus formulated as chance constraints with acceptable violation probability $\epsilon$. 
Since $\epsilon > 0$, there will be cases in which HVDC and PSTs reach their limit and are not able to continue to provide corrective control according to \eqref{HVDCcontrol}, \eqref{PSTcontrol}. This control saturation can be expressed as 
\begin{align}
    \tilde{p}_{DC,i} &=  \min \left\{ p_{DC,i} + \delta_{DC,i}^{kl} - \boldsymbol{\alpha}_{DC(i,\cdot)}\omega , ~~\bar{p}_{DC,i}\right\} , \\
    \tilde{\gamma}_{j} &=  \min \left\{ \gamma_{j} + \delta_{\gamma,j}^{kl} - \boldsymbol{\alpha}_{\gamma(j,\cdot)}\omega , ~~\bar{\gamma}_{j}\right\} .
\end{align}
When such saturation occurs, we get different power flows than expected, and there may be overloading on the transmission lines, requiring further corrective action. However, our chance constraints ensure that the need for such additional corrective action occurs with probability smaller than $\epsilon$.

\subsubsection{Power Flow Constraints}
The power flow constraints can be expressed as 
\begin{align}
& \mathbb{P}\left[ p_{ij}^{\delta_{kl}} +\! LF_{ij}^{kl} p_{kl}^{\delta_{kl}} +\! \left(\mathbf{A}_{(ij,\cdot)}\! + LF_{ij}^{kl}\mathbf{A}_{(kl,\cdot)}\right) \omega \leq \bar{p}_{ij}\right]\! \geq\! 1\!-\!\epsilon_l,  \nonumber\\
& \mathbb{P}\left[ p_{ij}^{\delta_{kl}} +\! LF_{ij}^{kl} p_{kl}^{\delta_{kl}} +\! \left(\mathbf{A}_{(ij,\cdot)}\! + LF_{ij}^{kl}\mathbf{A}_{(kl,\cdot)}\right) \omega \geq \underline{p}_{ij} \right]\! \geq\! 1\!-\!\epsilon_l, \nonumber\\
& \quad\quad\quad\quad\quad\quad\quad\quad\quad\quad\quad\quad\quad \forall_{ij \in \mathcal{L}, ~ kl \in \{0,\mathcal{L}\}} \label{chance_line_min}
\end{align}
where $\bar{p}_{ij}, \underline{p}_{ij}$ are the upper and lower limits on the power flow, enforced as chance constraints with acceptable violation probability $\epsilon_l$. Index $kl=0$ refers to the pre-contingency constraint with $LF_{ij}^{kl}=0$.  
Different from the limits on generators, HVDC and PSTs, the transmission limits are soft constraints from a physical perspective and might experience violations. The violation probability $\epsilon_l$ is thus interpreted as the probability of observing an overload or N-1 violation. These such transmission line violations can either be accepted (e.g., if the overload is small or related to an unlikely outage), or can be handled through additional measures in real-time.  

Note that the problem does not include any consideration of the intermediate post-contingency, pre-corrective control system state, as described in e.g. \cite{capitanescu2007}. Constraints to ensure feasibility of this state can however be added without any conceptual change to the method. 

\subsection{Reformulation of Chance Constraints}
The chance constraints in Eq.~\eqref{chance_gen_min}-\eqref{chance_line_min} need to be reformulated into tractable constraints. Here, we choose to use an analytic reformulation based on the assumption that the forecast uncertainty vector $\omega$ follows a normal distribution, similar to the one used in \cite{roald2013,  bienstock2014}. 
This choice of distribution is based on the fact that the chance constraints used are marginal chance constraints and are one dimensional (e.g., one constraint for each line flow) and depend on the weighted sum of the high dimensional random vector $\omega$. With the projection of a high dimensional vector onto a one dimensional constraint, arguments similar to the central limit theorem can be invoked \cite{dasgupta2006}, despite the fact that the individual uncertainty sources $\omega$ may not be normally distributed.

\subsubsection{Reformulated Constraints} 
With the normality assumption, Eqs.~\eqref{chance_gen_min}-\eqref{chance_line_min} can be reformulated into the following constraints: 
\begin{align}
&\alpha_{G,i}\Phi^{-1}(1-\epsilon_G)\sigma_{\Omega} \leq r_i^{+}, \label{lin_gen_max}\\
&\alpha_{G,i}\Phi^{-1}(1-\epsilon_G)\sigma_{\Omega} \geq -r_i^{-}, \label{lin_gen_min}
\\
&p_{DC,i} + \delta_{DC,i}^{kl} + \Phi^{-1}(1-\epsilon)\parallel\Sigma_W^{1/2} \boldsymbol{\alpha}_{DC(i,\cdot)}^T \parallel_2 \leq \bar{p}_{DC,i} \label{SOC_dc_max}\\
&p_{DC,i} + \delta_{DC,i}^{kl} - \Phi^{-1}(1-\epsilon)\parallel\Sigma_W^{1/2} \boldsymbol{\alpha}_{DC(i,\cdot)}^T \parallel_2 \geq \underline{p}_{DC,i} \label{SOC_dc_min} \\
&\gamma_{j} + \delta_{\gamma,j}^{kl} + \Phi^{-1}(1-\epsilon)\parallel \Sigma_W^{1/2} \boldsymbol{\alpha}_{\gamma(j,\cdot)}^{T}\parallel_2 \leq \bar{\gamma}_{j} \label{SOC_pst_max}\\
&\gamma_{j} + \delta_{\gamma,j}^{kl} - \Phi^{-1}(1-\epsilon)\parallel \Sigma_W^{1/2} \boldsymbol{\alpha}_{\gamma(j,\cdot)}^{T}\parallel_2 \geq \underline{\gamma}_{j}\label{SOC_pst_min} \\
& p_{ij}^{\delta_{kl}} + LF_{ij}^{kl} p_{kl}^{\delta_{kl}} ~+ \label{SOC_line_max}\\
&\quad\quad\Phi^{-1}(1-\epsilon_l) \parallel \Sigma_W^{1/2} \left(\mathbf{A}_{(ij,\cdot)} + LF_{ij}^{kl}\mathbf{A}_{(kl,\cdot)}\right)^{T} \parallel_2 \leq \bar{p}_{ij}  \nonumber\\
& p_{ij}^{\delta_{kl}} + LF_{ij}^{kl} p_{kl}^{\delta_{kl}} ~- \label{SOC_line_min}\\
&\quad\quad\Phi^{-1}(1-\epsilon_l) \parallel \Sigma_W^{1/2} \left(\mathbf{A}_{(ij,\cdot)} + LF_{ij}^{kl}\mathbf{A}_{(kl,\cdot)}\right)^{T} \parallel_2 \geq -\bar{p}_{ij} \nonumber
\end{align}
The generator constraints are linear since they only depend on the total power mismatch $\Omega$ with standard deviation $\sigma_{\Omega}$ \cite{li2015}. The remaining constraints \eqref{SOC_dc_max} - \eqref{SOC_line_min} are Second Order Cone (SOC) constraints, which are convex for $\epsilon\leq0.5$ \cite{bienstock2014}.

\subsubsection{Exploiting Symmetry of SOCs}
Each pair of upper and lower SOC constraints in \eqref{SOC_dc_max}-\eqref{SOC_line_min} can be reduced to a pair of linear constraints and a single SOC constraint by exploiting the symmetry of the normal distribution \cite{bienstock2014}:
\begin{align}
& p_{DC} + \delta_{DC}^{kl} + s_{DC} \leq p_{DC}^{max} \label{final_DC_max}\\
& p_{DC} + \delta_{DC}^{kl} - s_{DC} \geq p_{DC}^{min} \label{final_DC_min}\\
& s_{DC,i} \geq \Phi^{-1}(1-\epsilon)\parallel\Sigma_W^{1/2} \boldsymbol{\alpha}_{DC(i,\cdot)}^T \parallel_2, ~~\forall_{i\in \mathcal{H}} \label{final_DC_SOC} 
\end{align}
\begin{align}
& \gamma + \delta_{\gamma}^{kl} + s_{\gamma} \leq \gamma^{max} \label{final_PST_max} \\
& \gamma + \delta_{\gamma}^{kl} - s_{\gamma} \geq \gamma^{min} \label{final_PST_min} \\
& s_{\gamma,j} \geq \Phi^{-1}(1-\epsilon)\parallel \Sigma_W^{1/2} \boldsymbol{\alpha}_{\gamma(j,\cdot)}^{T}\parallel_2, ~~~\forall_{j\in \mathcal{S}} \label{final_PST_SOC}\\
&p_{ij}^{\delta_{kl}} + LF_{ij}^{kl} p_{kl}^{\delta_{kl}} + s_{ij}^{kl} \leq p_{ij}^{max}, \quad\quad~~~ \forall_{ij \in \mathcal{L}, ~ kl \in \{0,\mathcal{L}\}} \label{final_line_max}\\
&p_{ij}^{\delta_{kl}} + LF_{ij}^{kl} p_{kl}^{\delta_{kl}} - s_{ij}^{kl} \geq -p_{ij}^{max}, \quad~~~~ \forall_{ij \in \mathcal{L}, ~ kl \in \{0,\mathcal{L}\}}\label{final_line_min}\\
& s_{ij}^{kl} \geq \Phi^{-1}(1-\epsilon_l) \parallel \Sigma_W^{1/2} \left(\mathbf{A}_{(ij,\cdot)} + LF_{ij}^{kl}\mathbf{A}_{(kl,\cdot)}\right)^{T} \parallel_2 \label{final_line_SOC}
\end{align}
The above reformulation cuts the number of SOC constraints in half and thus improves efficiency. 
Note that the SOC terms \eqref{final_DC_SOC}, \eqref{final_PST_SOC}, \eqref{final_line_SOC} are always positive, and introduce a tightening of the original, deterministic constraints.

\subsubsection{Non-Normal Uncertainty} For cases where the assumption of a normal distribution is not justified,  an analytical reformulation can be still be obtained even when only limited knowledge of the distribution is available \cite{roaldArxiv, summers2014}. For example, a reformulation based on the Chebyshev inequality only requires knowledge of the mean and covariance of $\omega$. These reformulations result in SOC constraints similar to \eqref{final_DC_max}-\eqref{final_line_SOC}, 
and can be easily used within the suggested framework. However, the distributionally robust reformulations can be overly conservative \cite{roaldArxiv}, leading to unnecessarily high cost or infeasibility.
Instead, we suggest to use the normal assumption (as described above) combined with out-of-sample testing.



 \section{Solution Algorithm}
\label{sec:implementation}

The full OPF problem with security and chance constraints in Section \ref{sec:CCOPF} is a Second Order Cone Program (SOCP), with the SOC constraints given by Eq.~\eqref{final_DC_SOC}, \eqref{final_PST_SOC} and \eqref{final_line_SOC}. 
The problem has a large number of linear and SOC constraints, due to consideration of both the pre- and post-contingency situations. 
Although the SOC constraints are convex, it has been observed in the literature \cite{bienstock2014} that attempting to solve the entire optimization problem at once by using a non-linear SOCP solver can result in unacceptably long convergence times and run into numerical difficulties. The state of the art method in the literature for solving such chance constrained OPF problems is the sequential outer approximation cutting planes algorithm \cite{bienstock2014}. 
In this algorithm, a relaxed version of the problem is solved by eliminating all the SOC constraints. Then the relaxation is successively tightened using a sequence of polyhedral outer approximations obtained by adding separating hyperplanes (cutting-planes) to the violated SOC constraints until the desired accuracy is reached. The success of the algorithm relies on exploiting three salient properties prevalent in the CC-OPF \cite{bienstock2014}: (i) only a small fraction of the non-linear chance constraints (SOCs) are active at optimality, and these correspond to critical/congested lines, (ii) each SOC term depends only on a limited number of decision variables, and (iii) linear programming solvers have better speed and stability compared to non-linear solvers. 

However, there are some critical differences in the features of the CC-SCOPF considered in this paper that makes the cutting-plane algorithm unsuitable: 

\emph{Feature 1:} Due to the SOC constraints for the HVDC and PSTs, that are often tight at optimality, as well as the security constraints for the lines, more constraints must be approximated through polyhedral constraints. 

\emph{Feature 2:} Each SOC term, in particular the SOC terms for the lines, depend on a large number of decision variables since $\alpha_G,~\boldsymbol{\alpha}_{DC},~\boldsymbol{\alpha}_{\gamma}$ are potentially large matrices. 
For those constraints, we found that the cutting-plane algorithm requires a large number of iterations to obtain a good polyhedral approximation and a feasible solution within acceptable tolerance.
This implies that the optimization problem needs to be solved a many times, and that a significant amount of memory is required to store the large and increasing number of linear constraints.

\emph{Feature 3:} Due to security constraints, there is a very large number of SOCs present in the problem. Evaluating the SOC constraints at any given candidate solution is therefore time consuming. Since the cutting-planes algorithm requires a large number of iterations with an SOC evaluation in each iteration, the problem solves very slowly.

To overcome these difficulties and obtain an efficient implementation, we develop a sequential SOCP algorithm based on a constraint generation process customized to our problem.
As in \cite{bienstock2014}, we first solve a relaxed version of the problem involving only linear constraints. Instead of adding cutting-planes, we then add the full SOC terms for the most violated constraints. We continue solving a sequence of SOCPs until all constraints are satisfied, reaching the globally optimal solution of the original problem. The sequence in which the violated constraints are added has a significant impact on the solving time, and needs to be carefully chosen. 
In the following we describe the details of the algorithm and the reasoning behind. \\
\emph{\underline{Preprocessing:} Solving the SCOPF Without Uncertainty}\\
As a pre-processing step, we solve the SCOPF without consideration of uncertainty. This allows us to obtain a fast, first estimate of the active constraints.\\
\emph{\underline{Step 1:} Solving the CC-OPF Without Security Constraints}
\begin{itemize}
\item[(a)] We first solve a base case problem consisting of the power balance and generator constraints \eqref{nodalbal}-\eqref{res_cont}, \eqref{lin_gen_max}, \eqref{lin_gen_min}, the full PST and HVDC constraints \eqref{final_DC_max} - \eqref{final_PST_SOC} as well as the linear pre-contingency line constraints \eqref{final_line_max}, \eqref{final_line_min}. Since most of the SOC constraints for HVDC and PST are tight and there are few of them (i.e., no additional SOC terms for the security constraints), adding the full SOC upfront eliminates unnecessary iterations.
In addition, the SOC terms belonging to the pre-contingency constraints violated in the SCOPF are added.

\item[(b)] The line SOC constraints for the base case, i.e., \eqref{final_line_SOC} with $kl = 0$ are then checked for violations, and the most violated ones are added to the problem with a warm start from the previous iteration. This process is repeated until all the base case constraints are satisfied.
\end{itemize}
\emph{\underline{Step 2:} Solving the Full CC-SCOPF With Security Constraints}\\
We check for violation of the security constraints and add them sequentially using warm start. This part of the algorithm runs in three phases: 
\begin{itemize}
\item[(c)] In the first iteration, we add all post-contingency constraints that were active in the SCOPF without uncertainty, as these most likely will be active in the CC-SCOPF as well.
We include both the linear constraints \eqref{final_line_max}, \eqref{final_line_min} and the corresponding SOC constraint \eqref{final_line_SOC} to the problem. 
\item[(d)] After the first iteration, we check for violation of only the linear security constraints \eqref{final_line_max}, \eqref{final_line_min}, which is much faster than evaluating the full SOC constraints. 
However, since the SOC constraint implies a tightening of the linear constraint, a violation of the linear constraint almost always means that the corresponding SOC is violated as well. 
For the most violated constraints, we therefore add both the linear constraints \eqref{final_line_max}, \eqref{final_line_min} and the corresponding SOC constraint \eqref{final_line_SOC} to the problem. 
This process is repeated until all the linear security constraints are satisfied.
\item[(e)] We then check for the violated SOC terms for the line security constraints \eqref{final_line_SOC}, and add the most violated ones to the problem. As mentioned in \emph{Feature 3}, the number of post-contingency SOCs are rather large, and hence inefficient to evaluate. However, we observe that this stage of the algorithm only requires few iterations until all constraints are satisfied, since most violated security constraints were added in (d).\\
To reduce the computational time involved in each iteration (d), we do a pre-screening of the SOC terms based on the LODF matrix. Even though most $LF_{ij}^{kl}$ are non-zero, most are very small $<1e-3$. 
By evaluating the post-contingency SOC constraints only for those pairs $ij$ and $kl$ where $LF_{ij}^{kl}$ exceeds a certain threshold, the number of evaluations can be significantly reduced while maintaining acceptable accuracy. 
\end{itemize}
Finally, we check whether the current solution violates any of the pre-contingency constraints. If yes, we restart from (b). Otherwise, the algorithm is terminated and a globally optimal solution which satisfies all the constraints of the full problem has been found. In all encountered cases, only one pass of the algorithm was required to find an optimal solution without requiring to return to (b).

\section{Case study - IEEE 118 Bus System}
\label{sec:Case}

In this section, we analyze the benefits of corrective control based on a case study for the IEEE 118 bus system. The scalability of the proposed method and the sequential SOCP algorithm is demonstrated in the next section. 

\subsection{IEEE 118 Bus Test System}
\begin{figure}
\includegraphics[width=0.8\columnwidth]{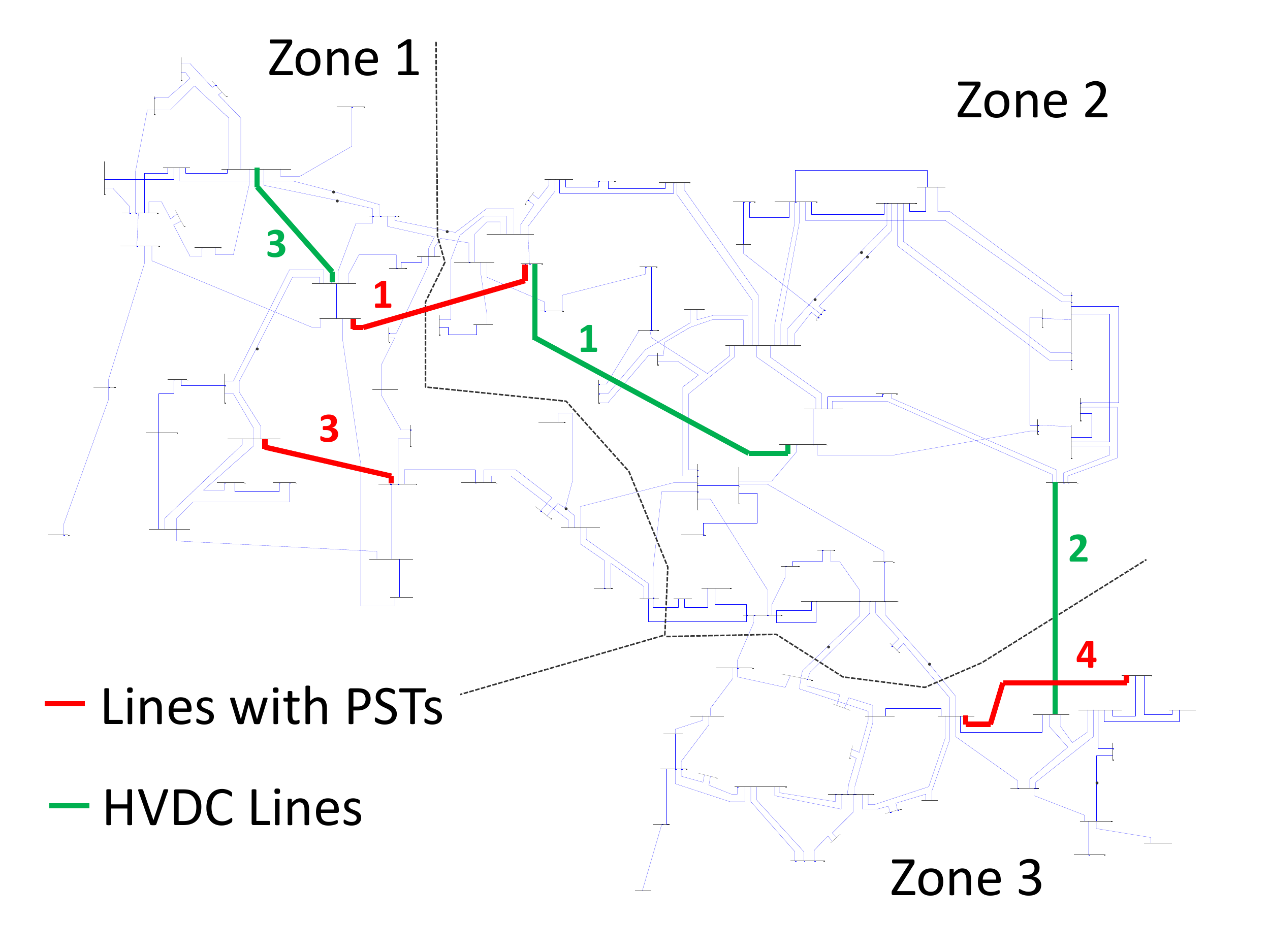}
\centering
\caption{IEEE 118 bus system with 3 Zones. The the lines with PSTs are marked in red, while HVDC connections are drawn in green.}
\label{system}
\end{figure}

We use the IEEE 118 bus test system as defined in \cite{118busdata}, with the following modifications. Both load and maximum generation capacity are scaled by a factor of 1.25, and the minimum generation capacity is set to zero.
The system loads are interpreted as a mix between load and renewable energy sources connected at a lower voltage level. Instead of adding uncertain in-feeds to the system, we assume that all 99 loads fluctuate around their forecasted consumption. The standard deviation $\sigma$ of each load is equal to 10\% of the forecasted consumption. As shown in Fig. \ref{system}, the system is divided into three zones. We assume that fluctuations within a zone are correlated with $\rho=0.3$, that fluctuations in different zones are uncorrelated. 
For the chance constraints, we apply $\epsilon_l=\epsilon=0.01$ for transmission line, HVDC and PST constraints. 
For the generator constraints, we use $\epsilon_G=0.001$, which corresponds to the acceptable probability of having insufficient reserves in the Swiss power system. 

We assume that there are 3 PSTs and 3 HVDC connections installed in the system, as shown in Fig. \ref{system}. Details for each HVDC connection are listed in Table \ref{table_HVDC_PST}. The PSTs are installed on lines 41, 167 and 54, with maximum tap positions $\pm 30^{\circ}$. The upper bounds on the post-contingency control for HVDC and PSTs are set to $\overline{\delta}_{DC}^{ij}=0.25~\overline{p}_{DC}$ and $\overline{\delta}_{\gamma}^{ij}=0.25~\overline{\gamma}$.

{
\begin{table} 
\caption{IEEE 118 Bus System - HVDC Connections}
\label{table_HVDC_PST}
\centering
\begin{tabular}{|l|l|l|l|}
\hline
 HVDC connection    & HVDC 1    & HVDC 2    & HVDC 3\\
\hline
From - To Bus       & 38 - 65   & 104 - 62  & 12 - 17\\
Replaced Line       & 96        & -         & 20\\
Capacity [MW]       & 500       & 200       & 175\\
\hline
\end{tabular}
\end{table}}

The total required down-reserves is set to $R^-=\Phi^{-1}(1-\epsilon_G)\sigma_\Omega$, which ensures that the total fluctuation $\Omega$ will be covered with probability $1-\epsilon_G$. The total amount of up-reserves $R^+$ is required to cover the same fluctuation or the maximum generation outage, whichever is larger: $R^+=\max\{\overline{p}_G, R^-\}$. This reserve dimensioning ensures that a similar amount of reserves are procured for the chance-constrained and deterministic formulations we will compare. 
To ensure ramping constraints, the upper bound on provision of reserves from each generator is set to 20\% of total generation capacity, i.e. $\overline{r}^+=\overline{r}^-=0.2\overline{p}_G$.

To analyze the performance of the method, we use a Monte Carlo simulation. We sample the uncertainty realization $\omega$, calculate the power flows and record the number of constraint violations. To test the algorithm for different types of uncertainty, we run two simulations with 2000 samples of $\omega$ each:

First, we use uncertainty data which correspond to our assumptions, and draw samples from a multivariate normal distribution with mean $\mu=0$ and covariance matrix $\Sigma_W$.

Second, we run out-of-sample tests based on 1 year of real system data from the Austrian Power Grid (APG).
The deviations are defined based on the difference between the the so-called DACF (Day-Ahead Congestion Forecast) and the snapshot (the real-time power injections) for all hours and buses with available data (8492 data points for 28 buses). Splitting the data into four three-month periods, we obtain 2000 data samples for each of the 99 buses. The samples are then scaled to fit the assumed covariance matrix $\Sigma_W$ and mean $\mu=0$. This corresponds to the case where we have a good estimate of the mean and covariance, but do not know the full underlying distribution.




We implement the sequential SOCP algorithm in Julia with JuMP \cite{jump}, and solve the problem using MOSEK \cite{mosek}. 
With this set-up, the solution is obtained within a few minutes on a laptop. 

\subsection{Numerical Results}
To demonstrate the benefits of corrective control, we compare five different solutions:
\begin{itemize}
    \item[(a)] The \emph{OPF} does not consider uncertainty (i.e., $\Sigma_W=0$) or constraints related to outages.
    \item[(b)] The \emph{SCOPF} does not consider uncertainty (i.e., $\Sigma_W=0$) or post-contingency corrective control ($\delta_{DC}^{ij}=\delta_{\gamma}^{ij}=0$).
    \item[(c)] The \emph{post-contingency corrective SCOPF} does not consider uncertainty (i.e., $\Sigma_W=0$).
    \item[(d)] The \emph{post-contingency corrective CC-SCOPF} does not consider corrective control in reaction to fluctuations (i.e., $\alpha_{DC}=\alpha_{\gamma}=0$).
    \item[(e)] The \emph{full CC-SCOPF} considers corrective control for both fluctuations and contingencies.
\end{itemize}
We first compare the cost of the generation dispatch and reserves (as obtained directly from the optimization), and analyze the differences. 
Second, we look into the empirical number of violations observed in the two Monte Carlo simulations for the normally distributed and out-of-sample data. 
Finally, we investigate the impact of corrective control by comparing the solution of (d) and (e) for different acceptable violation probabilities $\epsilon$ and varying uncertainty levels $\sigma$.




 

\subsubsection{Comparison of operation cost}
\begin{figure}
\includegraphics[width=0.98\columnwidth]{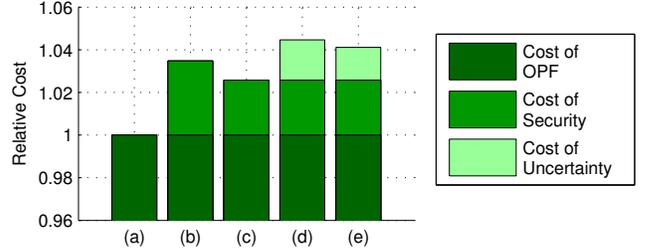}
\centering
\caption{Breakdown of the cost for the five OPF formulations (a) - (e). The costs are normalized by the cost of the OPF (a).}
\label{cost}
\end{figure}

The costs of the porblems (a)-(e) are listed in Table \ref{tableII} and are shown graphically in Fig. \ref{cost}, where the costs are normalized by the cost of the standard OPF.
The cost increase is analyzed with respect to the following criteria \cite{roald2013}:

\emph{Cost of Security} is defined as the cost increase between the OPF and the SCOPF, due to the cost of enforcing N-1 constraints.

\emph{Cost of Uncertainty} is defined as the cost increase between the SCOPF and the CC-SCOPF, due to the cost of enforcing chance constraints instead of deterministic constraints. 

We observe that the cost of security can be reduced by 0.9\% from (b) to (c), by introducing post-contingency corrective control. Similarly, the cost of uncertainty is reduced by 0.35\% from (d) to (e), by introducing uncertainty corrective control.

\begin{figure}
\includegraphics[width=0.98\columnwidth]{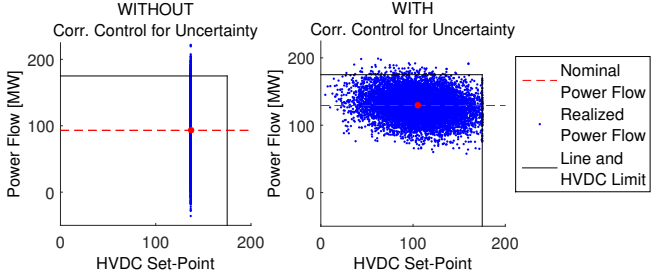}
\centering
\caption{Power Flow on Line 119 after outage of Line 126, plotted against the set-point of HVDC 3 for the CC-SCOPF without (left) and with (right) corrective control for uncertainty. The red marker shows the forecasted operating point, the blue points are the realized power flows obtained through the Monte Carlo simulation and the black lines show the line and HVDC limits.}
\label{lineflow}
\end{figure}

\begin{figure}
\includegraphics[width=0.95\columnwidth]{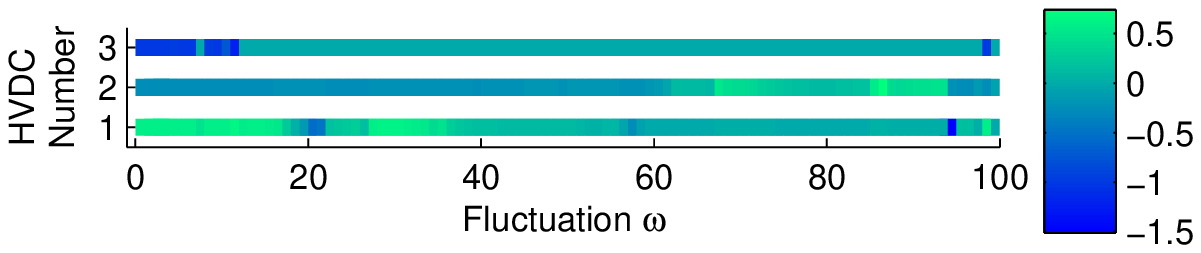}
\includegraphics[width=0.95\columnwidth]{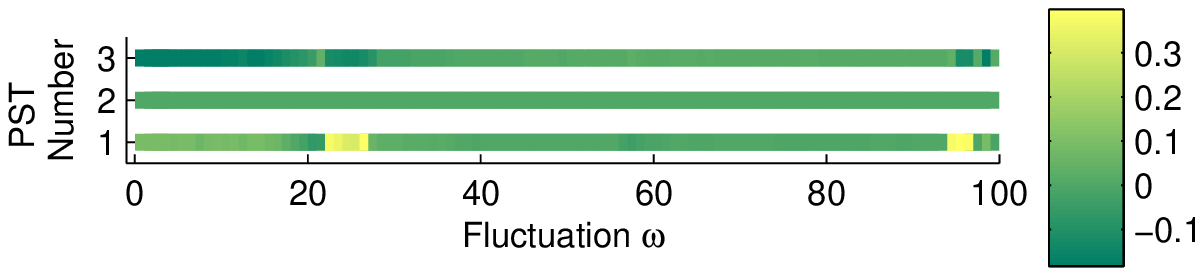}
\centering
\caption{Values of the corrective control parameters $\boldsymbol{\alpha}_{DC}$ (top) and $\boldsymbol{\alpha}_{\gamma}$ (bottom).}
\label{DC_PST}
\end{figure}

The reduction in cost of uncertainty can be explained by the better ability to react to the fluctuations in the power injections. 
In Fig. \ref{lineflow}, the power flow on line 119 after the outage of line 126 is plotted against the set-point of HVDC 3 for the CC-SCOPF problems without (left) and with (right) corrective control. The red marker shows the forecasted operating point, and the blue points correspond to the actual, realized conditions as simulated for 2000 samples of $\omega$. In the case without corrective control (left), the HVDC set-point remains constant for all uncertainty samples, while the set-point changes depending on the wind condition with corrective control (right).
The post-contingency constraint on line 119 is one of the active constraints in the CC-SCOPF problems (d) and (e), and it is clearly seen how the line flow limit is violated in some cases.
However, we also observe that by introducing the corrective control of HVDC and PSTs, it is possible to reduce the standard deviation of the line flow, in this case from 35.2 MW to 19.7 MW. With reduced variance, a higher nominal power flow can be accepted, which leads to a reduction in cost. 

The response of the HVDC and PSTs to uncertainty is determined by the corrective control parameters $\boldsymbol{\alpha}_{DC},~\boldsymbol{\alpha}_{\gamma}$. The values of those variables as obtained from (e) are shown in Fig. \ref{DC_PST}. We observe that the PST and HVDC react strongly to fluctuations that are geographically close to where they are located, and have close to zero entries for most other fluctuations.



\subsubsection{Number of Violations and Out-of-Sample Testing}
We run two Monte Carlo simulations based on the normally distributed and APG data sets, and compute the number of constraint violations. 

\begin{figure}
\includegraphics[width=1.0\columnwidth]{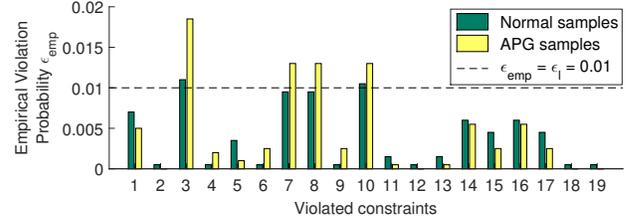}
\centering
\caption{Empirical violation probability $\epsilon_{emp}$ per constraint with normal samples (green) and APG samples (yellow) for the 19 violated constraints, based on the solution for the full CC-SCOPF (e). The black dotted line corresponds to $\epsilon_{emp}=\epsilon_l=0.01$.}
\label{ViolPerLine}
\end{figure}

To investigate how well our solutions perform in the out-of-samples test with real APG data, we first look at the empirical violation probability $\epsilon_{emp}$ for each separate constraint and compare this to the accepted violation probability $\epsilon_l$.
Fig. \ref{ViolPerLine} shows the empirical violation probability per constraint for the full CC-SCOPF (e), with normal (green bar) and APG (yellow bar) samples. Only the 19 constraints with $\epsilon_{emp}>0$ are shown.
We observe that there are four tight constraints (3, 7, 8 and 10) for which $\epsilon_{emp}\approx0.01$ with normally distributed samples. 
For these constriants, the APG samples lead to $\epsilon_{emp}>0.01$, indicating that the chance constraint is violated. However, the empirical violation probability remains below 2\% for all constraints, and below the acceptable 1\% value for the majority of constraints. While the empirical violation probability will vary depending on the uncertainty data, this result indicates that a chance constraint approach based on the normal distribution still significantly reduces the risk of violations.

In a second investigation, we calculate the joint violation probability, i.e., the probability that a sample will exhibit a violation of \emph{any} constraint. 
Since each chance constraint in our problem is enforced individually with a violation probability of $1-\epsilon$ and the wind deviations that lead to a constraint violation differ between constraints, the joint violation probability might be significantly higher than $1-\epsilon$.


{
\setlength{\extrarowheight}{2pt}
\begin{table} 
\caption{Cost and Number of Violations for (a) - (e)}
\label{tableII}
\centering
\begin{tabular}{l c c c c c c }
   & (a) & (b) & (c) & (d) & (e)\\[2pt]
\hline
Cost [\$]         &   84924    & 87880 & 87108 & 88713 & 88418\\ 
\hline
& & & & & \\[-2pt]
Violations & (a) & (b) & (c) & (d) & (e)\\[+2pt]
\hline
Normal      &  2000   &  1729    &  1712    &  98     &  87    \\[-2pt]
Samples  & (100\%) & (86.5\%) & (85.6\%) & (4.9\%) & (4.4\%) \\[+2pt]
\hline
APG        &  2000   &  1697    &  1655    &  114    &  95    \\[-2pt] 
Samples    & (100\%) & (84.9\%) & (82.8\%) & (5.7\%) & (4.8\%) \\[+2pt]
\hline
\end{tabular}
\end{table}}

In Table \ref{tableII}, the number of Monte Carlo samples which lead to either pre- or post-contingency constraint violations are listed for both the normally distributed and the APG samples.
The OPF (a) violates all scenarios, as it neither accounts for outages nor uncertainty. The two SCOPF formulations without (b) and with (c) corrective control have violations due to uncertainty in around 1700 of the samples. This corresponds to a joint violation probability of around 85\%, which is unacceptably high and clearly shows the need to account for uncertainty.
 Both CC-SCOPF formulations (d) and (e) have violations in around 100 samples, indicating a joint violation probability of around 5\%. As expected, this is above the acceptable violation probability $\epsilon_l=0.01$ for each separate constraint. However, it is still in the acceptable range and shows that a CC-SCOPF based on separate chance constraints also effectively reduces the joint violation probability.

The Monte Carlo simulation based on APG data increases number of violated scenarios by less than 1 \% compared to the normally distributed samples for the CC-SCOPF cases (a) and (e), indicating that the method performs well also in out-of-sample tests.
We further conclude that the CC-SCOPF with corrective control (e) is able to provide a similarly low violation probability as the CC-SCOPF without corrective control (d), but at lower cost.

\subsubsection{Influence of the confidence level} Table \ref{tableIII} shows the cost of the generation dispatch without (d) and with (e) corrective control for uncertainty, for different acceptable violation probabilities $\epsilon$ and different standard deviations (in \% of forecasted load). We observe that the possibility to react to uncertainty becomes increasingly important (i.e., the reduction in cost is larger) if we want to enforce smaller violation probabilities or if the level of uncertainty increases. With standard deviations of 14.25\%, we even observe that the CC-SCOPF with corrective control (e) is feasible, whereas the CC-SCOPF without (d) is not.

{
\setlength{\extrarowheight}{2pt}
\begin{table} 
\caption{Cost of CC-SCOPF without (d) and with (e) corrective control for uncertainty, for different values of $\epsilon$ (top) and different standard deviations $\sigma$ as \% of forecasted load (bottom).}
\label{tableIII}
\centering
\begin{tabular}{l c c c c c c }
      & $\epsilon\!=\!0.1$ & $\epsilon\!=\!0.05$ & $\epsilon\!=\!0.02$ & $\epsilon\!=\!0.01$ & $\epsilon\!=\!0.001$\\
\hline
(d)    & 87898 & 88136 & 88466 & 88713 & 89483\\
(e)       & 87720 & 87941 & 88201 & 88418 & 89070 \\
\hline
\\
            & $\sigma\!=\!5\%$ & $\sigma\!=\!7.5\%$ & $\sigma\!=\!10\%$  & $\sigma\!=\!12.5\%$ & $\sigma\!=\!14.25\%$\\
\hline
(d) &  84171       & 86334          & 88136          &  91513          & Infeasible\\
(e)    &  83991       & 86119          & 87941          &  91022          & 93063\\
\hline
\end{tabular}
\end{table}}

\section{Case study - Scalability}
\label{sec:Polish}
In this case study, we discuss the scalability of the CC-SCOPF and proposed solution algorithm based on the run times for three different test cases, the IEEE 118 bus system, the IEEE 300 bus system and the Polish test case with 2383 buses.

\subsection{Test System Data}
\subsubsection{IEEE 118 bus test system} For the IEEE 118 bus system, we use the same data as in the previous case study.
\subsubsection{IEEE 300 bus test system}
We use a modified version of the IEEE 300 bus test system found in the NICTA Energy System Test Archive \cite{coffrin2014}. To obtain a feasible SCOPF case, we increase the line limits by a factor of 5.0 and the generation limits by a factor of 2.0. 
All loads above 50 MW, corresponding to 109 loads and 68\% of the system load, are assumed to be uncertain. Their standard deviations are set to 5\% of the forecasted consumption, and the correlation to zero.
We assume that there are 3 HVDC connections installed in the system, with corresponding system data listed in Table \ref{table_HVDC_PST_300}. The PSTs are installed on lines 91, 140 and 174, with maximum tap positions $\pm 30^{\circ}$.


{
\begin{table} 
\caption{IEEE 300 Bus System - HVDC Connections}
\label{table_HVDC_PST_300}
\centering
\begin{tabular}{|l|l|l|l|}
\hline
 HVDC connection    & HVDC 1    & HVDC 2    & HVDC 3\\
\hline
From - To Bus       & 68 - 198  & 8 - 18    & 126 - 145\\
Replaced Line       & -         & -         & -\\
Capacity [MW]       & 900       & 600       & 600\\
\hline
\end{tabular}
\end{table}}

\subsubsection{Polish Test System}
We base our case study on the Polish Winter Peak test case with 2383 buses as provided with Matpower 5.1 \cite{zimmermann2011}. To obtain a feasible case with active transmission constraints, the maximum generation capacity is scaled by a factor of 2.0 and the transmission line limits by a factor of 2.5.
All loads above 25 MW, corresponding to 157 loads and 29\% of the system load, are assumed to be uncertain. Their standard deviations are set to 20\% of the forecasted consumption, and the correlation between loads to zero.
We assume that there are 3 HVDC connections installed in the system, with corresponding system data listed in Table \ref{table_HVDC_PST_Polish}. The PSTs are installed on lines 130, 240 and 1381, with maximum tap positions $\pm 30^{\circ}$.


{
\begin{table} 
\caption{Polish Test Case - HVDC Connections}
\label{table_HVDC_PST_Polish}
\centering
\begin{tabular}{|l|l|l|l|}
\hline
 HVDC connection    & HVDC 1    & HVDC 2    & HVDC 3\\
\hline
From - To Bus       & 32 - 18   & 184 - 105 & 67 - 138\\
Replaced Line       & -         & -         & 169\\
Capacity [MW]       & 500       & 500       & 1000\\
\hline
\end{tabular}
\end{table}}

For both of the above test systems, we choose similar settings as for the IEEE 118 bus test case in the previous section. We apply $\epsilon_l=\epsilon=0.01$ for transmission line, HVDC and PST constraints and use $\epsilon_G=0.001$ for generation constraints. We determine the total required reserve capacities $R^-,~R^+$ and the reserve provision bounds $\overline{r}^+,~\overline{r}^-$ as described above. We do not include post-contingency corrective control in either of the test cases, but solve the CC-SCOPF with corrective control for uncertainties. 
The sequential SOCP algorithm is implemented in Julia with JuMP \cite{jump}, and solved 
using Ipopt \cite{ipopt}. 

\subsection{Results}
The size of each test case is listed in Table \ref{table_size}.
Due to the large number of security constraints, already the IEEE 118 bus test case has more SOC constraints than the Polish test case solved in \cite{bienstock2014}. The problem sizes for the IEEE 300 bus and Polish grid test cases are one and two orders of magnitude larger than the IEEE 118 bus system, as measured in number of constraints, with the Polish test case featuring more than 8 million variables and 8 million SOC constraints.


Table \ref{table_times} shows the run times of the algorithm based on the current implementation. We have included both the total solution time, as well as an overview of the number of SOC evaluations and the total time spent evaluating the SOC constraints. 
We observe that a significant part of the solution time is spent checking the SOC constraints. The pre-screening of the security constraints based on the LODF matrix reduces the number of security constraints that we need to check by a factor of 3. However, for the largest Polish test case, we still need to evaluate close to 1.4 million SOC constraints. 

By using the sequential SOCP algorithm, we are able to keep the number of SOC evaluations very low, with only 2-3 evaluations. 
While the cutting plane algorithm \cite{bienstock2014} would solve the problem faster and more reliably in each iteration, it would require a much larger number of SOC constraint evalutions, leading to prohibitively large run times. However, the time required for the SOC check can be significantly reduced by more efficient coding and parallelization, and is left as a topic for future work. This would improve the solution times for both our solution algorithm and the cutting plane algorithm in \cite{bienstock2014}.

{
\begin{table} 
\caption{Size of the CC-SCOPF for the different test cases}
\label{table_size}
\centering
\begin{tabular}{l c c c}
                    & IEEE 118  & IEEE 300  & Polish\\[+1pt]
\hline
 &   &   &  \\[-8pt]
Variables           & ~54 000   & ~173 000  & $>$8 mill.\\
Linear Constraints  & ~70 000   & ~343 000  & $>$16 mill.\\
SOC Constraints     & ~34 000   & ~169 000  & $>$8 mill.\\[+1pt]
\hline
\end{tabular}
\end{table}}

{
\begin{table} 
\caption{Run times for CC-SCOPF for the different test cases}
\label{table_times}
\centering
\begin{tabular}{l l l l}
                            & IEEE 118  & IEEE 300   & Polish\\[+1pt]
\hline
 &   &   &  \\[-8pt]
Total Run Time              & 2min 15s  & 4 min 44s  & 4h 20min\\
Number of SOC Evaluations   & 2         & 3          & 2\\
Time per SOC Evaluation     & 8s         & 27 s       & 1h 36min\\
\hline
\end{tabular}
\end{table}}

\section{Conclusions}
\label{sec:Conclusion}

In this paper, we have proposed a framework for corrective control which involves corrective actions in response to both contingencies and different uncertainty realizations. \\
The combined corrective control was incorporated in a CC-SCOPF with security and chance constraints, to ensure that the probability of constraint violations does not exceed a desired level. The chance constraints were reformulated using an analytical approach, which leads to an optimization problem with SOC constraints. To be able to handle the problem computationally, we developed a sequential SOCP algorithm.\\
With this solution algorithm, we are able to solve the IEEE 118 bus case including contingency constraints and 99 uncertain loads within a few minutes on a laptop. The solutions obtained for the IEEE 118 bus system show that the use of corrective control actions in reaction to uncertainty reduces operational cost, without reducing the level of security in the system. 
The cost reduction is more significant for cases where we want to achieve a low violation probability, or where the level of uncertainty is higher. 

We also demonstrate the scalability of the method by presenting results for the IEEE 300 bus test system and the Polish grid with 2383 buses. One main bottleneck of the current implementation is the evaluation of the SOC terms, but we believe a better implementation and parallelization will allow for significant speed up of this part.  

In the case study, we observed that the PSTs and HVDC react similarly to fluctuations that are located close to each other, implying that geographically close fluctuations could be aggregated.
With this type aggregations, coupled with the fact that even large transmission systems only have few active contraints, we believe that the proposed CC-SCOPF is scalable to even larger systems with hundreds of uncertain variables. In addition to further scalability tests and improvements on the current implementation, future work will explore the extension to a full AC power flow, with approximate chance constraints based on a linearization of the power flow around the forecasted operating point \cite{qu15}.

\section{Acknowledgements}
This research work described in this paper has been partially carried out within the scope of the project "Innovative tools for future coordinated and stable operation of the pan-European electricity transmission system (UMBRELLA)", supported under the 7th Framework Programme of the European Union, grant agreement 282775. We thank our partners for useful discussions and feedback, and Austrian Power Grid for providing historical data.

\bibliographystyle{IEEEtran}
\bibliography{201512_bib.bib}

\vspace*{-2\baselineskip}

\begin{IEEEbiography}[{\includegraphics[width=1in,height=1.25in,clip,keepaspectratio]{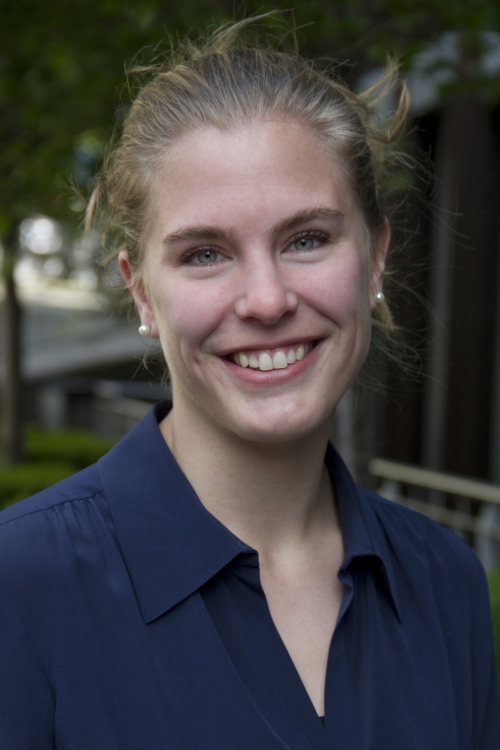}}]
{Line Roald}
(M'12) obtained her  B.Sc ('09) and M.Sc. ('12) in Mechanical Engineering from ETH Z\"urich, Switzerland. Currently, she is pursuing a PhD degree in the Power Systems Laboratory, Department of Electrical and Computer Engineering at ETH Z\"urich. Her research interests include modeling and optimization of power system operation under consideration of uncertainty and risk-based security criteria.
\end{IEEEbiography}

\vspace*{-2\baselineskip}

\begin{IEEEbiography}[{\includegraphics[width=1in,height=1.25in,clip,keepaspectratio]{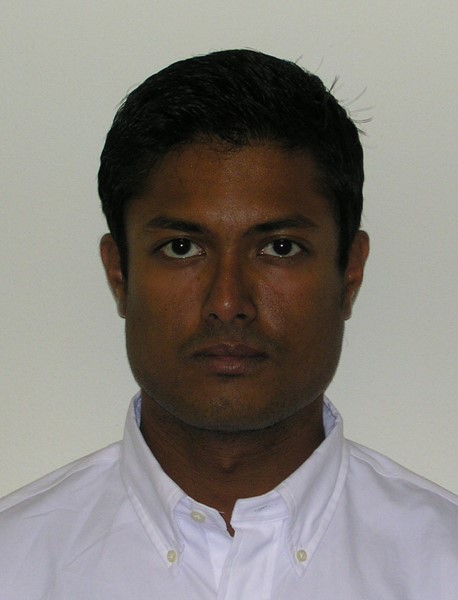}}]
{Sidhant Misra}
is a research scientist in the theory division at Los Alamos National Laboratory. He obtained his PhD degree in Electrical Engineering and Computer Science from Massachusetts Institute of Technology in 2014. His research interests include optimization, machine learning and random graph theory, with applications to risk based optimization and inference in energy networks.
\end{IEEEbiography}

\vspace*{-2\baselineskip}

\begin{IEEEbiography}[{\includegraphics[width=1in,height=1.25in,clip,keepaspectratio]{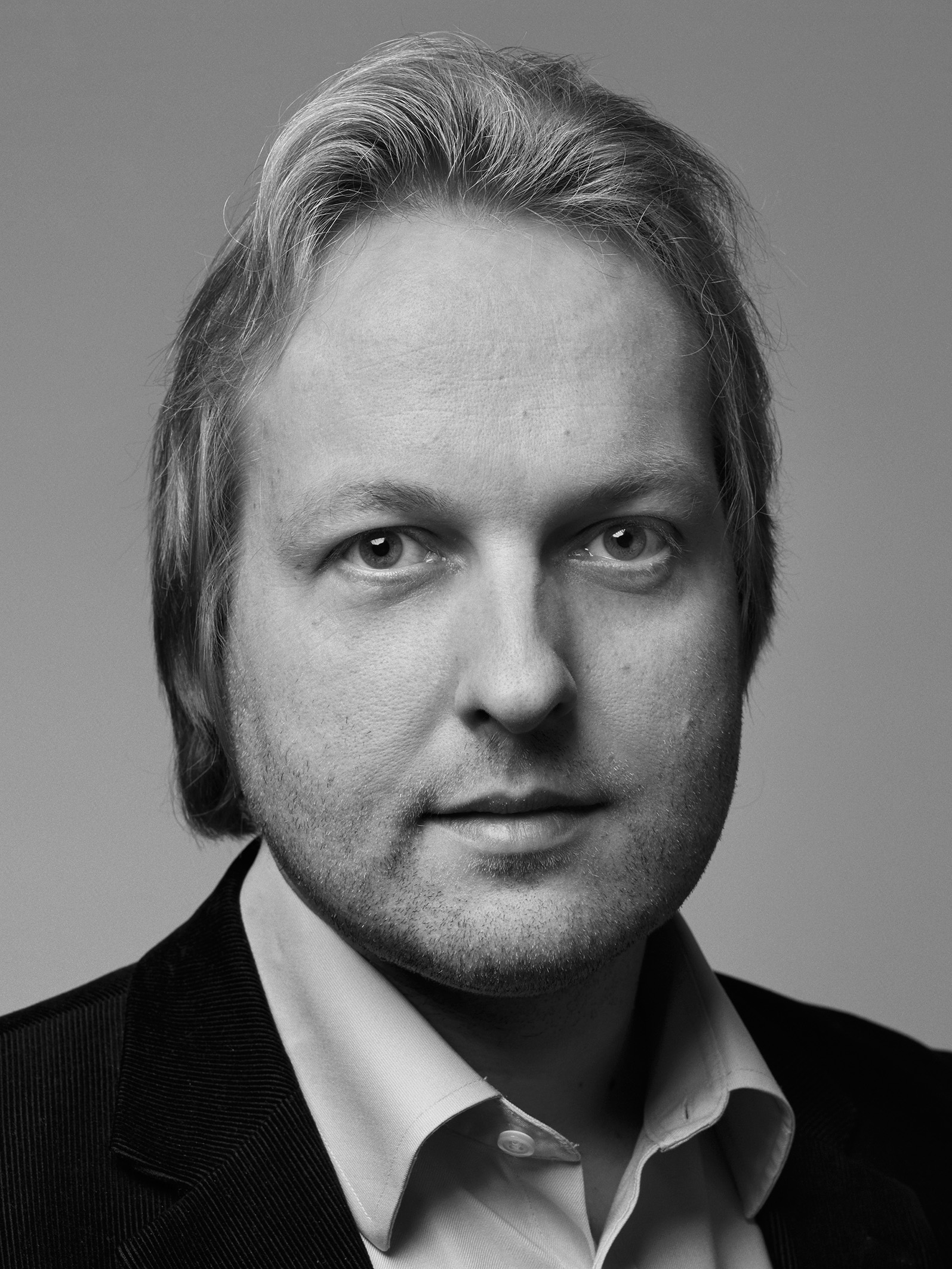}}]
{Thilo Krause}
received the Diplom-Wirtschafts-Ingenieur degree in economics and electrical engineering from the Dresden University of Technology, Dresden, Germany, in 2002 and the Ph.D. degree from the Swiss Federal Institute of Technology (ETH), Zürich, Switzerland, in 2006. Until 2016, he has been a senior researcher at the Power Systems Laboratory, ETH Zürich. Recently, he joined the strategic asset management section of the main electricity supplier of the City of Zurich (ewz).
\end{IEEEbiography}

\vspace*{-2\baselineskip}

\begin{IEEEbiography}[{\includegraphics[width=1in,height=1.25in,clip,keepaspectratio]{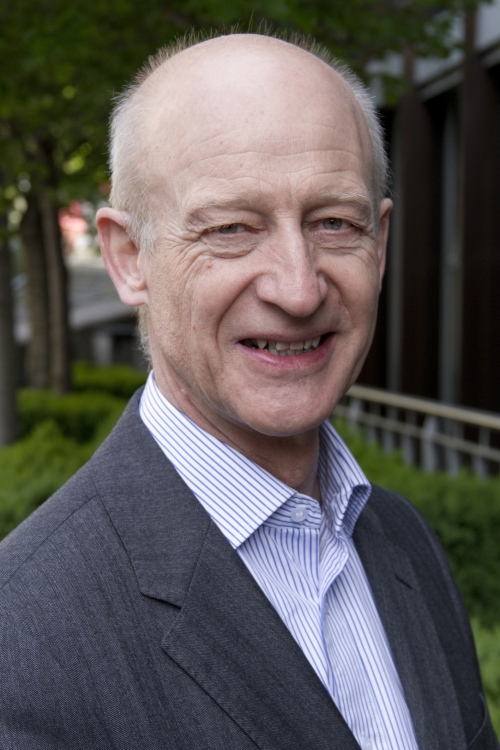}}]
{G\"oran Andersson}
(M’86, SM’91, F’97) obtained his M.S. (1975) and Ph.D. (1980) degrees from the University of Lund, Sweden. In 1980 he joined ASEA’s, now ABB’s, HVDC division in Ludvika, Sweden, and in 1986 he was appointed full professor in electric power systems at KTH (Royal Institute of Technology), Stockholm, Sweden. Since 2000 he is full professor in electric power systems at ETH Zürich (Swiss Federal Institute of Technology). His research interests include power system dynamics, control and operation, power markets, and future energy systems.
Göran Andersson is Fellow of the Royal Swedish Academy of Sciences, the Royal Swedish Academy of Engineering Sciences, the Swiss Academy of Engineering Sciences, and foreign member of the US National Academy of Engineering. He was the recipient of the 2007 IEEE PES Outstanding Power Educator Award, the 2010 George Montefiore International Award 2010, and the 2016 IEEE PES Prabha S. Kundur Power System Dynamics and Control Award.
\end{IEEEbiography}

\end{document}